\documentclass[11pt]{article}
\usepackage{amssymb}
\newcommand{\be}{\begin{equation}}
\newcommand{\ee}{\end{equation}}

\newcommand{\ep}{\epsilon}
\newtheorem{theorem}{Theorem}
\newtheorem{lemma}{Lemma}
\begin{document}
\title{ The KLR-theorem revisited}
\author{Abram Kagan \\
Department  of Mathematics, University of Maryland\\
 College Park, MD 20742, USA}
 \date{}
\maketitle
\begin{abstract}
\noindent
For independent random variables $X_1,\ldots, X_n;Y_1,\ldots, Y_n$  with all $X_i$ identically distributed and same for $Y_j$, we study the relation
\[E\{a\bar X + b\bar Y|X_1 -\bar X +Y_1 -\bar Y,\ldots,X_n -\bar X +Y_n -\bar Y\}={\rm const}\]  
with $a, b$ some constants.
It is proved that for $n\geq 3$ and $ab>0$ the relation holds iff $X_i$ and $Y_j$ are Gaussian.\\
A new characterization arises in case of $a=1, b= -1$. In this case either $X_i$ or $Y_j$ or both have a Gaussian component. It is the first (at least known to the author) case when presence of a Gaussian component is a characteristic property. 
\end{abstract}
\section{Introduction and an ancillary lemma}
Let $X_1,\ldots, X_n$ be independent identically distributed (iid) random variables with
$E(|X_i|)<\infty,\:\bar X=(X_1 +\ldots+X_n)/n$. The KLR-theorem (see Kagan {\it et al.}(1965), also Kagan {it et al.}(1973), Ch. 5)
claims that for $n\geq 3$ the relation
\be 
E\{\bar X|X_1 -\bar X,\ldots, X_n -\bar X\}={\rm const}
\ee
is characteristic of Gaussian $X_i$.\\

The following lemma is a version of the KLR-theorem.
\begin{lemma}Let $X_1,\ldots, X_n,Y_1,\ldots,Y_n$ be independent with $X_1,\ldots,X_n$
id with $E(|X_i|)<\infty$ and $Y_1,\ldots,Y_n$ id. Then for $n\geq 3$ the relation
\be
E\{\bar X|X_1-\bar X +Y_1-\bar Y,\ldots, X_n-\bar X +Y_n-\bar Y\}={\rm const}
\ee
holds iff $X_i$ is Gaussian.
\end{lemma}
{\it Proof of Lemma}. Replacing $X_i$ with $X_i-E(X_i)$ the constant in (2) may be made 
zero. On setting $f_{1}(t)=E\{\exp it X_j\},\:f_{2}(t)=E\{\exp it Y_j\}$,
 multiplying both parts of (2) by $\exp\sum_{1}^{n} t_{j}(X_j-\bar X+Y_{j}-\bar Y)$ and taking the expectation results in
\be 
\sum_{j=1}^{n}f'_{1}(t_j -\bar t)\prod_{k\neq j}f_{1}(t_k -\bar t)\prod_{j=1}^{n}f_{2}(t_j -\bar t)=0,
\ee
 where $\bar t =(t_1+\ldots+t_n)/n$.\\
 On setting $\tau_j=t_j -\bar t$ and dividing both side of (3) by $\prod_{1}^{n}f_{1}(\tau_j) f_{2}(\tau_j)$ for $|\tau_j|<\ep,\:\ep >0$ small enough where $f_{1}(\tau_j)f_{2}(\tau_j)\neq 0$ leads to
 \be
 \sum_{1}^{n}\frac{f'_1}{f'_2} (\tau_j)=0,\:\sum_{1}^{n}\tau_j =0,\:|\tau_j|<\ep.
 \ee
 This is the Cauchy classical functional equation whose solutions are linear functions,
 $f'_{1}(t)/f_{1}(t) =ct$ for some $c$ and $|t|<\ep$ whence 
 \be
 f_{1}(t)=\exp ct^2/2, |t|<\ep.
 \ee
 Since $f_{1}(t)$ is a characteristic function, $c<0$ and (5) holds for all $t$ proving the necessity part.\\
 As to the sufficiency part, if $X_i$ is Gaussian and $Y_i$ arbitrary, due to independence of $(\bar X, X_1-\bar X,\ldots, X_n -\bar X)$ and $(Y_1 -\bar Y,\ldots,Y_n-\bar Y)$,one has
 \be 
 E\{\bar X|X_1 -\bar X, Y_1 -\bar Y,\ldots, X_n -\bar X, Y_n -\bar Y\}=E\{\bar X|X_1-\bar X,\ldots, X_n -\bar X\}=0.
 \ee
 Since the $\sigma-$algebra $\sigma(X_1-\bar X+Y_1 -\bar Y,\ldots, X_n -\bar X +Y_n -\bar Y)$ is a subalgebra of $\sigma(X_1 -\bar X,Y_1-\bar Y,\ldots, X_n -\bar X, Y_n -\bar Y)$, one has due to (6)
 \be
 E\{\bar X|X_1 -\bar X +Y_1 -\bar Y,\ldots, X_n -\bar X+Y_n -\bar Y\}=0.
 \ee
 \section{The main result}
 Here we are studying the relation
 \be
 E\{a\bar X + b \bar Y|X_1 -\bar X +Y_1 -\bar Y,\ldots, X_n -\bar X +Y_n -\bar Y\}={\rm const}
 \ee
 with the $a, b$ some constants
 and show that when $ab>0$ and $n\geq 3$ it is characteristic for Gaussian $X_i$ and $Y_i$ while in case of $ab<0$ a new phenomenon arises.
 \begin{theorem} Let $X_1,\ldots, X_n, Y_1,\ldots, Y_n$ be as in Lemma with \\
 $E\{(X_i)^2\}<\infty,\\ E\{(Y_i)^2\}<\infty $.\\
 For $n\geq 3$ and $ab>0$ the relation (8) holds iff $X_i$ and $Y_i$ are Gaussian.
 \end{theorem}
 {\it Proof of Theorem 1}. As in Lemma, one may assume without loss in generality $E(X_i)=E(Y_i)=0$. Also one may always assume $a>0, b>0$.\\
 Proceeding as in the proof of Lemma, one gets from (8) a functional equation for the characteristic functions $f_{1}(t),\:f_{2}(t)$  of $X_i$ and $Y_i$:
 \be
 a\sum_{j=1}^{n}f'_{1}(\tau_j)\prod_{k\neq j}f_{1}(\tau_k)\prod_{1}^{n}f_{2}(\tau_j)+
 b\sum_{j=1}^{n}f'_{2}(\tau_j)\prod_{k\neq j}f_{2}(\tau_k)\prod_{1}^{n}f_{1}(\tau_j)=0
 \ee
 holding for the $\tau_1,\ldots,\tau_n$ with $\sum_{1}^{n}\tau_i =0,\:|\tau_i|<\ep$
 with $\ep >0$ as in Lemma.\\
 Dividing both sides of (9) by $\prod f_{1}(\tau_i)f_{2}(\tau_i)$ and setting 
  $g(\tau)=af'_{1}(\tau)/f_{1}(\tau) + bf'_{2}(\tau)/f_{2}(\tau)$ results in the Cauchy
  functional equation
  \be
  \sum_{1}^{n}g(\tau_i)=0\: {\rm if}\:\sum_{1}^{n}\tau_i =0, |\tau_i|<\ep
  \ee 
  whence for $g(t)=ct$ for some $c$ and $|t|<\ep$. From (10)
  \be 
  \{f_{1}(t)\}^a \{f_{2}(t)\}^b =\exp ct^{2}/2, |t|<\ep.
  \ee
  From $g'(0)=af''_{1}+bf''_{2}(0) = -aE\{(X_i)^2\} -bE\{(Y_i)^2\}=c$, one has $c<0$. 
  \begin{lemma} Let $f_{1}t),\ldots, f_{m}(t)$ be characteristic functions and
  $\alpha_1,\ldots,\alpha_m$ positive numbers. The relation
  \be 
  \prod_{1}^{m}\{f_{i}(t)\}^{\alpha_i}=\exp{ct^2}, |t|<\ep
   \ee
   holds for some $c<0,\ep>0$ iff $f_{1}(t),\ldots, f_{m}(t)$ are characteristic functions of Gaussian distributions.
   \end{lemma}
   For a proof of Lemma 2 generalizing Cramer classical theorem on the components of
   a Gaussian distribution see Zinger and Linnik (1955) (see also Kagan {\it et al.}(1973, Ch. 1). \\
   Applying Lemma 2 to (11) completes the proof of Theorem 1.\\
   \\
   Turn now to the case when in (8) $a=1, b= -1$. Before stating the result, remind that
   a random variable $X$ is said having a Gaussian component if $X$ is equidistributed
   with $U+\xi$, denoted $X\cong U + \xi$, where $U$ and $\xi$ are independent and $\xi$ is (non-degenerate) Gaussian.
   \begin{theorem} Let $X_1,\ldots, X_n,Y_1,\ldots, Y_n$ be as in Theorem 1. Assume that the characteristic functions $f_{1}(t),\:f_{2}(t)$ of $X_i$ and $Y_i$ do not vanish.
   \\
   The relation
   \begin{eqnarray}
   E\{\bar X|X_1 -\bar X +Y_1 -\bar Y,\ldots, X_n -\bar X +Y_n -\bar Y\}=\nonumber \\
   E\{\bar Y|X_1 -\bar X +Y_1 -\bar Y,\ldots, X_n -\bar X +Y_n -\bar Y\}
   \end{eqnarray} 
   holds for $n\geq 3$ iff either (i) $X_i$ and $Y_i$ have the same (arbitrary) distribution
   or (ii) $Y_i$ is arbitrary (with $f_{2}(t)\neq 0$) while $X_i\cong Y_i +\xi$
   with $Y_i$ and  Gaussian $\xi$ independent or (iii) $X_i$  is arbitrary while $Y_i\cong X_i +\eta$ with $X_i$ and Gaussian $\eta$ independent. In other words, unless $X_i$ and $Y_i$
   are identically distributed, at least one of them has a Gaussian component. 
 \end{theorem} 
   {\it Proof of Theorem 2}. Proceeding as in proof of Theorem 1 and taking into account the condition $\prod_{1}^{n}f_{1}(\tau_i)f_{2}(\tau_i)\neq 0$ one gets the following 
   Cauchy equation for $g(t)=f'_{1}(t)/f_{1}(t) -f'_{2}(t)/f_{2}$:
   \be 
   \sum_{1}^{n}g(\tau_i)=0 \:\: {\rm if}\:\:\sum_{1}^{n}\tau_i =0.
   \ee
   From (14),
   \be f_{1}(t)=f_{2}(t)\exp ct^{2}/2
   \ee
   with $c=f''_{1}(0) - f''_{2}(0)= E\{(Y_i)^2\}-E\{(X_i)^2\}$.\\
   If $c=0$, is $X_i$ and $Y_i$ are identically distributed.
    If $c<0$, $X_i$ has a Gaussian component $\xi$ with variance $c$, $X_i\cong Y_i +\xi$ while $Y_i$ is arbitrary. If $c>0$, $Y_i$ has a Gaussian component $\eta$ with variance $-c$, $Y_i\cong X_i +\eta$ while $X_i$ is arbitrary. 
   This proves Theorem 2.\\
   \\
   In Kagan and Klebanov (2010) were described some analytic properties characterizing random variables with Gaussian components. It would be interesting to find a statistical property characteristic for distributions with Gaussian components. Here is such an example.\\
   \\
   Set $s_{X}^2=(1/n-1)\sum_{1}^n (X_i -\bar X)^2,s_{Y}^2=(1/n-1)\sum_{1}^n (Y_i -\bar Y)^2$. In Lukacs (1942) was shown that if $n\geq 2$ the relation
   \be E(s_{X}^2|\bar X)={\rm const}
   \ee
   holds iff  $X_i$ is Gaussian. A version of (16) in spirit of (13) is 
   \be
   E(s_{X}^2 |\bar X +\bar Y)=E(s_{Y}^2|\bar X  +\bar Y)+c
   \ee
   where $c=var(X_i)- var (Y_i)$.
   \begin{theorem}
   Under the conditions of Theorem 2, the relation (17) holds for $n\geq 2$ if (i) in case of $c>0$ $X_i\cong Y_i +\xi$ where $\xi$ is independent of $Y_i$ Gaussian random variable with variance $c$, (ii) in case of $c<0$ $Y_i\cong X_i +\eta$ where $\eta$ is independent of $X_i$ Gaussian random variable with variance $-c$, (iii) in case of $c=0$
   $X_i\cong Y_i$ (remind that we assumed $E(X_i)=E(Y_i)=0$).
   \end{theorem}
   {\it Proof of Theorem 3}. Multiplying both sides of (17) by  $\exp \{it(\bar X +\bar Y)\}$, taking the expectation and then dividing both sides by $(f_{1}(t))^n (f_{2}(t))^n$
   results in a differential equation for the characteristic functions $f_{1}(t),\:f_{2}(t)$ of $X_i$ and $Y_i$:
   \begin{eqnarray} 
   [f''_{1}(t) (f_{1}(t))^{n-1}-(f'_{1}(t))^2 (f_{1}(t))^{n-2}]/(f_{1}(t))^2 \nonumber \\
   =[f''_{2}(t) (f_{2}(t))^{n-1}-(f'_{2}(t))^2 (f_{2}(t))^{n-2}]/(f_{2}(t))^2 -c
   \end{eqnarray}
   equivalent to 
   \be
   \{\log f_{1}(t)\}'' =\{\log f_{2}(t)\}'' -c
   \ee
   whence  $f_{1}(t)=f_{2}(t)\exp (-ct^2/2)$ proving Theorem 3.
   
   \section*{References}
   Lukacs, E. A characterization of he normal distribution. {\it Ann. Math. Stat.}, {\bf 13}(1942), 1, 91-93.\\
   \\
   Zinger, A. A. And Linnik, Yu. V. On an analytic generalization of the Cramer theorem
   in Russia). {\it Vestnik of Leningrad Univ.}, 11 (1955), 51-56.\\
   \\
   Kagan, A. M., Linnik, Yu. V., Rao, C. R. On a characterization of the normal law based on a property of the sample average. {\it Sankhya}, {\bf A 27}, 3-4 (1965), 405-406.\\
  \\
   Kagan, A. M., Linnik, Yu. V., Rao, C. R. {\it Cjaracterization Problems in Mathematical Statistics}, Wiley, N. Y. (1973).\\
   \\
   Kagan, A. M., Klebanov, L. B. A class of multivariate distributions related to distributions with a Gaussian component. In: {\it IMS collections}, {\bf 7}(2010), 105-112.
 \end{document}